\documentclass[11pt]{article}
\usepackage[a4paper,margin=23mm]{geometry}
\usepackage{amsmath,amssymb,amsthm,microtype}
\usepackage[colorlinks=true,
linkcolor=blue,citecolor=blue,
urlcolor=blue]{hyperref}
\usepackage{xurl}
\newtheorem{theorem}{Theorem}
\newtheorem{subtheorem}{Theorem}
\newtheorem{conjecture}{Conjecture}
\newtheorem{claim}{Claim}
\newtheorem{proposition}{Proposition}
\newtheorem{construction}{Construction}

\newtheorem{lemma}{Lemma}[section]
\begin{document}

\title{Removable forests in $k$-connected and $k$-edge-connected graphs} 
\author{
    Xu Liu\thanks{College of Computer Science, Nankai University, Tianjin 300350, P.~R. China.}\and
Bo Ning\thanks{Corresponding author. College of Computer Science, Nankai University, Tianjin 300350, P.~R. China. \nolinkurl{Email:bo.ning@nankai.edu.cn}. Partially supported by the National Natural Science Foundation of China (No. 12371350).~\href{https://orcid.org/0000-0002-9622-5567}{ORCID: 0000-0002-9622-5567}.}}
\maketitle
\begin{abstract}
We prove the forest conjecture of Clay and Jord\'an: Let $k\ge1$ be an integer and let $F$ be a forest with $m$ vertices and $e$ edges. If a $k$-connected graph $G$ has minimum degree at least $\max\{k+e,m-1\}$, then it contains a copy of $F$ whose edge deletion preserves $k$-connectivity. We also prove the $k$-edge-connectivity version.
\end{abstract}


\section{Introduction}
A famous conjecture due to Lov\'asz~\cite{lovasz} in 1975, states that for every positive integer $k$ there is an integer $f(k)$ with the following property: Given two distinct vertices $x$ and $y$ in an $f(k)$-connected graph $G$, one can find an $(x,y)$-path $P$ such that $G-V(P)$ is $k$-connected. Tutte's theorem~\cite{tutte} gives $f(1)=3$. Chen, Gould, and Yu~\cite{chen-gould-yu} and, independently, Kriesell~\cite{kriesell} showed that $f(2)=5$ suffices. Later, Kawarabayashi and Ozeki~\cite{ko} asked how much connectivity can be retained after deleting the vertices of several paths with the same ends. They proved that for every $s\ge1$ and any two distinct vertices $x,y$ of a $(2s+1)$-connected graph, there are $s$ internally disjoint $(x,y)$-paths such that deleting all their vertices leaves a connected graph. The bound $2s+1$ is sharp. If the graph is $(3s+2)$-connected, the paths can be chosen so that the remaining graph is $2$-connected. A different question is whether one can prescribe the order of the subgraph to be deleted. Fujita and Kawarabayashi~\cite[Theorem~2]{fujita-kawarabayashi} showed that every $k$-connected graph $G$ with $\delta(G)\ge\lfloor3k/2\rfloor+2$ has an edge $uv$ such that $G-\{u,v\}$ is $k$-connected. Confirming a conjecture of Fujita and Kawarabayashi~\cite{fujita-kawarabayashi},  Mader~\cite{mader-paths} proved that a $k$-connected graph with minimum degree at least $\lfloor3k/2\rfloor+m-1$ contains an $m$-vertex path whose vertex deletion leaves a $k$-connected graph. He conjectured that the same bound suffices for a copy of any given tree on $m$ vertices. 

\begin{conjecture}[Mader~\cite{mader-paths}]\label{Conj:1}
    Let $k$ and $m$ be positive integers, and let $T$ be a tree of order $m$.
If $G$ is a $k$-connected graph with $\delta(G)\ge \lfloor 3k/2\rfloor+m-1$, then $G$ contains a subtree $T'\cong T$ such that $G-V(T')$ is $k$-connected.
\end{conjecture}

The case $k=1$ of Conjecture \ref{Conj:1} follows from a theorem of Diwan and Tholiya~\cite{diwan-tholiya}. Hong and Liu~\cite{hong-liu} proved the cases $k=2$ and $k=3$. For general $k,m\ge2$, Mader~\cite{mader-trees} already proved that minimum degree $2(k+m-1)^2+m-1$ suffices. Liu, Ying, and Hong~\cite{liu-ying-hong} lowered this bound to a linear one of $3k+4m-6$.

There are also other problems with similar studying style.
Thomassen~\cite{thomassen} studied vertex deletion for cycles. He proved that every $(k+3)$-connected graph $G$ has a cycle $C$ for which $G-V(C)$ is $k$-connected. For trees in bipartite graphs, Luo, Tian, and Wu~\cite{luo-tian-wu} conjectured that if $T$ has bipartition $(X,Y)$, then every $k$-connected bipartite graph $G$ with $\delta(G)\ge k+\max\{|X|,|Y|\}$ contains a copy $T'$ of $T$ such that $G-V(T')$ is $k$-connected. For some partial results, see \cite{YangTian2023,YangTian2025}. 
For a brief survey on Mader's conjecture, Luo-Tian-Wu Conjecture, and related results, see Section 4 of \cite{TianMeng2026}.

Following a similar line of study, in this article, we focus on the problem of deleting an edge set while preserving the 
$k$-connectivity of the resulting graph. We call a subgraph $H$ removable if $G-E(H)$ retains the required connectivity. Hasunuma~\cite[Conjecture~1.5]{hasunuma} conjectured that every $k$-connected or $k$-edge-connected graph of minimum degree at least $k+m-1$ contains a removable copy of each tree on $m$ vertices. 

\begin{conjecture}[Hasunuma~\cite{hasunuma}]\label{conj:edge, k+m-1}
For any $k\geq 1$ and any tree $T$ of order $m$, every $k$‐connected (respectively, $k$‐edge‐connected) graph $G$ with $\delta(G)\geq k+m-1$ contains a subtree $T'\simeq T$ such that $G-E(T')$ is $k$‐connected 
(respectively, $k$‐edge‐connected). 
\end{conjecture}
The following construction shows that the minimum degree condition in conjecture \ref{conj:edge, k+m-1} is sharp, we then prove in Section \ref{sec:construction}.
\begin{construction}\label{cons:tree}
    Let $k\ge1$ and $m\ge2$ be integers with $d=k+m-2\ge2$, and $n\ge d+1$ be an integer with $nd$ even. Let $r=\lfloor \frac{d}{2}\rfloor$. The graph $G$ has vertex set $V=\{v_0,...,v_{n-1}\}$. With subscripts modulo $n$, vertex $v_i$ is adjacent to $v_{i-j}$ and $v_{i+j}$ for each integer $1\le j\le r$. If $d$ is odd, we add the edges $v_iv_{i+\frac{n}{2}}$ into $E(G)$ for each integer $0\le i<\frac{n}{2}$. Let $T=K_{1,m-1}$ be a star.
\end{construction}

Hasunuma~\cite{hasunuma} proved both versions of Conjecture \ref{conj:edge, k+m-1} for $k=1,2$.
Liu, Liu, and Hong~\cite{liu} and, independently, Yang and Tian~\cite{YangTian2024} settled the case $k=3$. Recently, Clay and Jord\'an~\cite[Theorem~1.1]{clay} proved both versions for all $k$. 

They also considered similar problems for forests~\cite[Theorem~1.2]{clay}. Let $F$ be a forest with $m$ vertices and $\omega\le4$ components. They showed that $\delta(G)\ge k+m-\omega$ and $k\ge\omega-1$ guarantee a removable copy of $F$ in every $k$-connected graph $G$. The same holds for $k$-edge-connectivity. In the same paper \cite{clay}, Clay and Jord\'an  proposed the following conjecture.
\begin{conjecture}\label{Conj:3}
(i)    Let $G$ be a $k$-connected graph, and let $F$ be a forest of order $m$ with $\omega$ connected components. If
$\delta(G)\ge k+m-\omega$ and $k\ge\omega-1$, then $G$ contains a
subforest $F'\cong F$ such that $G-E(F')$ is $k$-connected.\\
(ii) Let $G$ be a $k$-edge-connected graph, and let $F$ be a forest of order $m$ with $\omega$ connected components. If
$\delta(G)\ge k+m-\omega$ and $k\ge\omega-1$, then $G$ contains a
subforest $F'\cong F$ such that $G-E(F')$ is $k$-edge-connected.
\end{conjecture}
We also have a construction, which implies the minimum degree condition in Conjecture \ref{Conj:3} is sharp.
\begin{construction}\label{cons:forest}
    Let $k\ge1$ and $m>\omega\ge1$ be integers with $k\ge \omega-1$. Let $e=m-\omega$ and $d=k+e-1\ge2$. Let $n\ge\max\{m,d+1\}$ be an integer with $nd$ even, $r=\lfloor\frac{d}{2}\rfloor$. The graph $G$ has vertex set $\{v_0,...,v_{n-1}\}$. With subscripts modulo $n$, vertex $v_i$ is adjacent to $v_{i-j}$ and $v_{i+j}$ for each integer $1\le j\le r$. If $d$ is odd, we add edges $v_iv_{i+\frac{n}{2}}$ for each integer $0\le i<\frac{n}{2}$. Let $F=K_{1,e}\cup (\omega-1)K_1$.
\end{construction}
We prove conjecture \ref{Conj:3} by proving the following theorem..
\begin{theorem}\label{thm:main}
Let $k\ge1$ be an integer, let $G$ be a graph, and let $F$ be a forest with $m$ vertices and $\omega$ components. Suppose that $\delta(G)\ge k+m-\omega$ and $k\ge\omega-1$.
\begin{enumerate}
\renewcommand{\labelenumi}{\textup{(\alph{enumi})}}
\item If $G$ is $k$-connected, it contains a copy $F'$ of $F$ such that $G-E(F')$ is $k$-connected.
\item If $G$ is $k$-edge-connected, it contains a copy $F'$ of $F$ such that $G-E(F')$ is $k$-edge-connected.
\end{enumerate}
\end{theorem}

We show that the degree condition in Theorem \ref{thm:main} is sharp. Let $e=e(F)=m-\omega$. We will prove that $\delta(G)\ge\max\{k+e,m-1\}$ suffices in both parts of Theorem~\ref{thm:main}. If $k\ge\omega-1$, then $k+e\ge m-1$, and the degree bound reduces to $k+m-\omega$. If $k<\omega-1$, then $m-1$ is the larger term.
For $e\ge1$, consider a star with $e$ edges in $K_{k+e}$. This complete graph has minimum degree $k+e-1$. Deleting the edges of any copy of the star leaves its center with degree $k-1$. But this contradicts the fact that both $k$-connectivity and $k$-edge-connectivity of a graph imply the minimum degree is at least $k$. When $m-1>k+e$, the complete graph $K_{m-1}$ has minimum degree $m-2$ and is both $k$-connected and $k$-edge-connected. It also has too few vertices to contain an $m$-vertex forest.

The minimum degree conditions in Conjecture \ref{Conj:3} and Theorem \ref{thm:main} are sharp in general. For $e\ge1$, let $F=K_{1,e}\cup(\omega-1)K_1$ and use Construction \ref{cons:forest}, whose connectivity properties are proved in Section \ref{sec:construction}. For $d\ge2$, deleting the edges of any copy of $F$ in $G$ leaves the centre of its star component with degree $d-e=k-1$, thus the remaining graph is neither $k$-connected nor $k$-edge-connected.

All graphs are finite, simple, and undirected. We follow Bondy and Murty~\cite[Chapters~1 and~2]{bondy-murty} and write $v(G)$, $e(G)$, and $\delta(G)$ for the order, size, and minimum degree of a graph $G$, respectively. A graph is $k$-connected if it has at least $k+1$ vertices and deleting fewer than $k$ vertices leaves a connected graph. A graph is $k$-edge-connected if it has at least two vertices and deleting fewer than $k$ edges leaves a connected graph. Let $S\subseteq V(G)$. We denote by $G[S]$ the subgraph of $G$ induced by $S$. For any two vertices $x,y\in V(G)$, an $(x,y)$-path means a path with $x,y$ as two end-vertices.
We emphasize that copies of trees and forests need not be induced.

In Section \ref{sec:construction}, we verify the tightness of Conjecture \ref{conj:edge, k+m-1} and \ref{Conj:3}. In Section~\ref{sec:vertex}, we prove Conjecture \ref{Conj:3} (i). We use Erd\H{o}s--S\'os theorem, a theorem of Mader and a theorem of Clay and Jord\'an. In Section~\ref{sec:edge}, we confirm Conjecture \ref{Conj:3}(ii). Our main tool includes Gomory--Hu tree.

\section{The proofs}
\subsection{Sharpness of the minimum degree condition.}\label{sec:construction}
\begin{proposition}
    The graphs in Construction \ref{cons:tree} and \ref{cons:forest} are $d$-regular, $d$-connected and $d$-edge-connected.
\end{proposition}
\begin{proof}
Let $G$ be either of the graphs defined in construction \ref{cons:tree} and \ref{cons:forest}. Since $n\ge d+1$, $r=\lfloor\frac{d}{2}\rfloor$. If $d$ is even, then each vertex $v_i$ in $G$ has degree $d(v_i)=2r=d$. If $d$ is odd, then $N(v_i)=\{v_{i-r},...,v_{i-1},v_{i+1},...,v_{i+r},v_{i+\frac{n}{2}}\}$, we have $d(v_i)=2r+1=d$. Thus $G$ is $d$-regular. 

Now we show that $G$ is $d$-connected (respectively, $d$-edge-connected). Suppose that $S\subseteq V(G)$ with $|S|\le d-1$. If $G-S$ is disconnected, let $V_1,\ldots,V_a$ be the maximal intervals of consecutive vertices of the cycle that remain in $G-S$, listed in cyclic order. For each $1\le i\le a$, $G[V_i]$ is connected. Assume that $V_i=\{v_{x_i},...,v_{y_i}\}$ and $V_{i+1}=\{v_{x_{i+1}},...,v_{y_{i+1}}\}$, where the vertices of $V_i$ have consecutive subscripts modulo $n$. Let $S_i=\{v_{y_i+1},...,v_{x_{i+1}-1}\}$. If there exist no edges between $V_i$ and $V_{i+1}$ in $G-S$, then $|S_i|\ge r$. 

Suppose that $d$ is even. If $V_i$ and $V_{i+1}$ lie in different components of $G-S$, then there exist at least two integers $i_1,i_2$, such that $|S_{i_1}|\ge r$ and $|S_{i_2}|\ge r$, since each $V_i$ is connected in $G-S$. Thus $|S|\ge|S_{i_1}|+|S_{i_2}|\ge2r$. However, $|S|\le d-1=2r-1$, a contradiction. 

Now suppose that $d=2r+1$. As above, there are at least two
distinct indices $i_1,i_2$ such that $|S_{i_1}|\ge r$
and $|S_{i_2}|\ge r$. Since $|S|\le d-1=2r$, we have $|S|=2r$, which implies $a=2$, and $V(G)-S=V_1\cup V_2$. Assume that $p=|V_1|\le|V_2|$ with $V_1=\{v_0,...,v_{p-1}\}$, $V_2=\{v_{p+r},...,v_{n-r-1}\}$. Since $p+r\le \frac{n}{2}$ and $\frac{n}{2}+p-1\le n-r-1$, for each $0\le i\le p-1$, we have $v_{i+\frac{n}{2}}\in V_2$ and $v_iv_{i+\frac{n}{2}}\in E(G-S)$,
which implies $V_1$ and $V_2$ are in the same component. Since $V_1$ and $V_2$ are both connected in $G-S$, we know that $G-S$ is connected, a contradiction.

Thus $G$ is $d$-connected. Since each $d$-connected graph is $d$-edge-connected, we complete the proof. 
\end{proof}
\subsection{The vertex-connectivity version}\label{sec:vertex}
For disjoint vertex sets $X$ and $Y$ of a graph $G$, let $e_{G}(X,Y)$ be the number of edges between them. Let $d_{G}(v)$ and $N_{G}(v)$ denote the degree and the set of neighbors of $v$ in $G$. For $v\in V(G)$ and $S\subseteq V(G)$, let $d_G(v,S)=|N_G(v)\cap S|$. When we apply Lemma~\ref{lem:insert}, $X$ is the set of vertices covered by the copies of the components already embedded.

\begin{subtheorem}\label{thm:vertex}
Let $G$ be $k$-connected, where $k\ge1$, and let $F$ be a forest with $m$ vertices and $e$ edges. If $\delta(G)\ge\max\{k+e,m-1\}$, then $G$ contains a copy $F'$ of $F$ for which $G-E(F')$ is $k$-connected.
\end{subtheorem}

We first state the three results used in the proof. The Erd\H{o}s--S\'os theorem will also be used in Section~\ref{sec:edge}.

\begin{lemma}[Clay and Jord\'an {\cite[Theorem~1.1]{clay}}]\label{lem:tree}
Let $k,t\ge1$ be integers, let $G$ be a $k$-connected graph with $\delta(G)\ge k+t$, and let $T$ be a tree with $t$ edges. Then $G$ contains a copy $T'$ of $T$ such that $G-E(T')$ is $k$-connected.
\end{lemma}

The original theorem is stated for a tree of order $m$. Lemma~\ref{lem:tree} follows by letting $m=t+1$. An edge of a $k$-connected graph is called \emph{essential} if deleting it destroys $k$-connectivity of the graph.

\begin{lemma}[Mader {\cite{mader}}]\label{lem:essential-cycle}
Let $G$ be a $k$-connected graph, where $k\ge1$. Every cycle of $G$ whose edges are all essential contains a vertex $v$ with $d_{G}(v)=k$.
\end{lemma}

Our proof needs the solution to Erd\H{o}s-S\'os in full by GPT-6
Astra \cite{GPT}. For two simplified proofs, see Riordan and Scott \cite{RiordanScott-arxiv-26}, and Wood {\cite[Theorem~1]{wood}}.
\begin{lemma}[Erd\H{o}s--S\'os]\label{lem:erdos-sos}
Let $s\ge2$ be an integer. Every nonempty graph of average degree greater than $s-2$ contains a copy of every tree on $s$ vertices.
\end{lemma}

Let $s=t+1$. By Lemma~\ref{lem:erdos-sos}, a graph $G$ containing no copy of a specified $t$-edge tree satisfies $2e(G)\le(t-1)v(G)$.

We also record the elementary criterion used to show that an edge is not essential.

\begin{lemma}\label{lem:replacement}
Let $G$ be a $k$-connected graph, where $k\ge1$, and let $uv\in E(G)$. If $G-uv$ contains $k$ internally vertex-disjoint $(u,v)$-paths, then $G-uv$ is $k$-connected.
\end{lemma}

\begin{proof}
Let $S\subseteq V(G)$ have size less than $k$. If $S$ contains $u$ or $v$, then $(G-uv)-S=G-S$ is connected. Otherwise at least one of the given paths contains no vertex of $S$. The graph $G-S$ is connected, and every occurrence of $uv$ in a path in $G-S$ can be replaced by this surviving $u,v$-path. Thus $(G-uv)-S$ is connected. Since $G-uv$ has the same vertex set as $G$, it is $k$-connected.
\end{proof}

\begin{lemma}\label{lem:insert}
Let $G$ be $k$-connected, where $k\ge1$, let $X\subsetneq V(G)$ have size $x$, and let $U=V(G)\setminus X$. Suppose that $d_G(v)\ge d$ for every $v\in U$. A given tree $T$ with $t\ge1$ edges has a copy in $G[U]$ whose edge deletion preserves $k$-connectivity under either of the following conditions.
\begin{enumerate}
\renewcommand{\labelenumi}{\textup{(\roman{enumi})}}
\item $d\ge x+t+1$, $d-k\ge2t$, and $2(d-k)\ge x+2t$.
\item $d=x+t$ and $x\ge2k$.
\end{enumerate}
\end{lemma}

\begin{proof}
Choose a $k$-connected spanning subgraph $H$ of $G$ with as few edges as possible, which contains each edge with an end in $X$. Let $H'=G[U]-E(H[U])$. If $H'$ contains a copy $T'$ of $T$, then $G-E(T')$ contains $H$ as a spanning subgraph, therefore $G-E(T')$ is $k$-connected. Suppose that $H'$ contains no copy of $T$. Let $Z$ consist of the vertices of $U$ of degree $k$ in $H$, and let $Y=U\setminus Z$, $y=|Y|$, $z=|Z|$. Every edge of $H[U]$ is essential in $H$, otherwise we could delete it and still keep all edges incident with $X$, contrary to the choice of $H$. Each vertex of $Y$ has degree greater than $k$ in $H$, then Lemma~\ref{lem:essential-cycle} implies that $H[Y]$ is a forest. Let $c$ be the number of components in $H[Y]$, with $c=0$ if $Y$ is empty.

We first sum $d_H(v)$ over $v\in Z$, which gives $kz$. When we sum $d_H(v)$ over $v\in Y$, the edges between $Y$ and $X$ contribute $e_H(X,Y)$. The edges of $H[Y]$ contribute $2(y-c)$, the edges between $Y$ and $Z$ contribute at most $kz$, since each has an end in $Z$, and each vertex of $Z$ has degree $k$ in $H$. Thus $\sum_{v\in U}d_H(v)\le e_H(X,Y)+2(y-c)+2kz$. 

Each edge of $G$ outside $H$ has both ends in $U$. Since $d_G(v)\ge d$ for each $v\in U$, we obtain
\begin{equation}\label{eq: 2e(Gamma)>=}
    2e(H')=\sum_{v\in U}(d_G(v)-d_H(v))\ge d|U|-e_H(X,Y)-2(y-c)-2kz.
\end{equation}

\noindent\textbf{Case 1.} (i) holds.

Each vertex of $Z$ has degree at least $d-k$ in $H'$, then $\sum_{v\in Z}d_{H'}(v)\ge(d-k)z$. On the other hand, the degree sum $\sum_{v\in Z}d_{H'}(v)$ is at most $e(H')+e(H'[Z])$, since an edge is counted twice only if both ends belong to $Z$. Thus $(d-k)z\le e(H')+e(H'[Z])$.

Applying Lemma~\ref{lem:erdos-sos} to $H'$ and $H'[Z]$, we have $2e(H')\le(t-1)(y+z)$ and $2e(H'[Z])\le(t-1)z$, thus 
$$
2(d-k)z\le 2(e(H')+e(H'[Z]))\le(t-1)(y+2z).
$$
If $t=1$, this inequality forces $Z=\varnothing$. If $t\ge2$, we have $2(d-k-t+1)z\le(t-1)y$. Since $d-k\ge2t$, we know $(d-k-t+1)\ge t+1>t-1$, thus $z\le y$. In both cases $Y$ is nonempty and $y\ge z$, thus $c\ge1$.

Since $|X|=x$, $e_H(X,Y)\le xy$. Since $2e(H')\le(t-1)|U|$, by \eqref{eq: 2e(Gamma)>=}, we have 
$$
d|U|-xy-2(y-c)-2kz\le 2e(H')\le (t-1)|U|.
$$
By $|U|=y+z$, we have $(d-x-t-1)y+(d-2k-t+1)z+2c\le0$. Let $y=(y-z)+z$, we obtain 
$(d-x-t-1)(y-z)+\bigl(2(d-k)-x-2t\bigr)z+2c\le0$. However, $2c>0$, and the first two terms are nonnegative by assumptions in (i), a contradiction. This completes the proof of Case 1.

\noindent\textbf{Case 2.} (ii) holds.

If $Y$ is empty, then $H'$ has minimum degree at least $d-k\ge t$, thus $H'$ contains a copy $T'$ of $T$. Hence we may assume $c\ge1$. We show that $e_H(X,Y)\le(x-1)y+c$.

Let $Y_0$ consist of the vertices of $Y$ adjacent to all of $X$ in $H$, and let $y_0=|Y_0|$. The set $Y_0$ is independent in $H$. Indeed, any edge inside it would have $k$ internally disjoint replacement paths through $X$, contradicting essentiality by Lemma~\ref{lem:replacement}.

\begin{claim}\label{Claim: d XcupY0 v}
    For each $v\in Y\setminus Y_0$, we have $d_H(v,X\cup Y_0)\le x-1$.
\end{claim}

\noindent\textbf{Proof of Claim \ref{Claim: d XcupY0 v}.}
This holds if $v$ has no neighbour in $Y_0$ immediately. Assume that $u\in N_H(v)\cap Y_0$, $p=\lvert N_H(v)\cap X\rvert$. If $p\ge k$, then the $p$ common neighbours in $X$ give $p\ge k$ internally disjoint paths of length two from $v$ to $u$. By Lemma~\ref{lem:replacement}, $H-vu$ is $k$-connected, contrary to the essentiality of $vu$. Thus $p<k$. 

If $v$ has more than $k$ neighbours in $X\cup Y_0$, then it has at least $k-p$ neighbours in $Y_0\setminus\{u\}$. Choose distinct vertices $a_1,\ldots,a_{k-p}$ in $X\setminus N_H(v)$ and distinct vertices $b_1,\ldots,b_{k-p}$ in $N_H(v)\cap(Y_0\setminus\{u\})$. For each $i$, the vertices $v,b_i,a_i,u$ form a path in $H-vu$, since both $b_i$ and $u$ are adjacent to all vertices in $X$. The new paths have distinct internal vertices, and their internal vertices in $X$ avoid those of the length-two paths $va'u$ where $a'\in N_H(v)\cap X$. This gives $k$ internally disjoint paths avoiding $vu$. By Lemma~\ref{lem:replacement}, $H-vu$ is $k$-connected, which  contradicts the essentiality of $vu$. Thus, $v$ has at most $k$ neighbours in $X\cup Y_0$. Since $x\ge 2k$ and $k\ge1$, we have $d_H(v,X\cup Y_0)\le k\le x-1$, the claim follows. $\hfill\blacksquare$

Since $H[Y]$ has $c$ components and $Y_0$ is a stable set, at least $y_0-c$ vertices of $Y_0$ have a neighbour in $Y\setminus Y_0$, thus $e_H(Y\setminus Y_0,Y_0)\ge y_0-c$. By Claim \ref{Claim: d XcupY0 v}, 
$$
\sum_{v\in Y\setminus Y_0}d_H(v,X\cup Y_0)=e_H(Y\setminus Y_0,X)+e_H(Y\setminus Y_0,Y_0)\le (x-1)(y-y_0).
$$
Then $e_H(Y\setminus Y_0,X)\le (x-1)(y-y_0)-y_0+c$. By the definition of $Y_0$, we have  
\begin{align*}
    e_H(X,Y)&\le(x-1)(y-y_0)-y_0+c+xy_0\\
    &=(x-1)y+c
\end{align*}
By $d=x+t$, $x\ge 2k$, $c\ge1$ and \eqref{eq: 2e(Gamma)>=}, we have
\begin{align*}
    2e(H')&\ge d|U|-e_H(X,Y)-2(y-c)-2kz\\
    &\ge (x+t)(y+z)-(x-1)y-c-2y+2c-2kz\\
    &=(t-1)y+(x+t-2k)z+c\\
    &\ge(t-1)(y+z)+z+c\\
    &>(t-1)|U|
\end{align*}
By Lemma~\ref{lem:erdos-sos}, $H'$ contains a copy of $T$, a contradiction. This completes the proof of Case 2, and hence of Lemma \ref{lem:insert}.
\end{proof}

\begin{proof}[Proof of Theorem~\ref{thm:vertex}]
If $e=0$, then $F=mK_1$, $G$ contains a copy of $F$. Suppose that $e\ge1$. Let $d=\max\{k+e,m-1\}$. We first embed all nontrivial components of $F$ and place the isolated vertices at the end. First embed a component of $F$ with the largest number of edges. If this component has $t$ edges, then $\delta(G)\ge d\ge k+e\ge k+t$. Lemma~\ref{lem:tree} gives a copy whose edge deletion preserves $k$-connectivity.

At each later step, we apply Lemma \ref{lem:insert} to the graph obtained from $G$ by deleting the edges of the previously chosen copies of tree components, and this graph is $k$-connected. Let $X$ be the vertex set of the $q$ components already embedded, and let $x=\lvert X\rvert$. Those components have $x-q$ edges. Since the $q$ components already embedded are nontrivial, we have $x\ge2q$. Let $t$ be the number of edges in the remaining largest component $T$. Each vertex outside $X$ still has degree at least $d$, since all previously deleted edges have both ends in $X$. 

The $q$ components already embedded have $x-q$ edges. Since $T$ is an unembedded component of $F$ and $x\ge2q$, we have $e\ge x-q+t\ge x/2+t$. Since $T$ has $t+1$ vertices and $|X|=x$, we have $m\ge x+t+1$, then $d\ge m-1\ge x+t$.

Suppose that $d>x+t$. Since the first component embedded also has at least $t$ edges and $d\ge k+e$, we have $d-k\ge e\ge2t$. Moreover, the $q$ components already embedded have $x-q$ edges in total, thus $2(d-k)\ge2e\ge x+2t$. Thus condition (i) of Lemma~\ref{lem:insert} holds. 

If $d=x+t$, then $x+t=d\ge k+e\ge k+x-q+t$, thus $q\ge k$ and $x\ge2k$. Thus condition (ii) holds. In either case we insert the next component on unused vertices and delete its edges while preserving $k$-connectivity. Repeat for all nontrivial components and then assign distinct unused vertices to the isolated vertices of $F$.
\end{proof}

\subsection{The edge-connectivity version}\label{sec:edge}
Throughout this subsection, $G$ is the host graph, and $F$ is the forest of order $m$ and edge count $e$. For $A\subseteq V(G)$, we denote by $\partial_G(A)$ the edge cut with sides $A$ and $V(G)\setminus A$, and $d_G(A)$ denotes its size.

\begin{subtheorem}\label{thm:edge}
Let $G$ be $k$-edge-connected, where $k\ge1$, and let $F$ be a forest with $m$ vertices and $e$ edges. If $\delta(G)\ge\max\{k+e,m-1\}$, then $G$ contains a copy $F'$ of $F$ for which $G-E(F')$ is $k$-edge-connected.
\end{subtheorem}

We use the Erd\H{o}s--S\'os theorem from Lemma~\ref{lem:erdos-sos} and the following form of the Gomory--Hu theorem~\cite{gomory}; see also Abboud, Krauthgamer, and Trabelsi~\cite[Theorem~1.1]{abboud-gh2021}. The edge weights are capacities and are assumed to be nonnegative, as in Gomory and Hu~\cite[p.~551]{gomory}. For a graph with nonnegative edge weights, the capacity of a cut is the sum of the weights of its edges. A minimum $(x,y)$-cut is a cut of least capacity separating $x$ and $y$.

\begin{lemma}[{\cite{gomory}; see also \cite{abboud-gh2021}}]\label{lem:cut-equivalent}
Every undirected graph $G$ with nonnegative edge weights
admits a tree $T$ on $V(G)$ with nonnegative edge weights
such that, for every pair of distinct vertices $x,y$,
each minimum $(x,y)$-cut in $T$ determines a minimum
$(x,y)$-cut in $G$ with the same capacity.
\end{lemma}

Here, the corresponding cuts in $T$ and $G$ are determined by the same bipartition of their common vertex set. Such a tree is called a \emph{cut-equivalent tree}.  Taking all edge weights of $G$ to be $1$ gives the following lemma.

\begin{lemma}\label{lem:gomory-hu}
Let $G$ be a connected graph with at least two vertices. There is a tree $T$ with vertex set $V(G)$ and positive integer edge weights such that the following properties hold.
\begin{enumerate}
\renewcommand{\labelenumi}{\textup{(\roman{enumi})}}
\item For distinct $x,y\in V(G)$, the minimum weight on the unique $(x,y)$-path in $T$ equals the minimum size of an edge cut separating $x$ and $y$ in $G$.
\item For every edge $xy$ of $T$, let $A$ be the vertex set of either component of $T-xy$. Then $d_G(A)$ equals the weight of $xy$.
\end{enumerate}
\end{lemma}

\begin{proof}
Give every edge of $G$ weight $1$ and choose a cut-equivalent tree $T$ by the theorem above. Write $w(e)$ for the weight of a tree edge $e$.

For distinct vertices $x,y$, deleting a minimum-weight
edge on their path in $T$ gives a minimum $(x,y)$-cut
in $T$. The same vertex partition gives a minimum $(x,y)$-cut in $G$ of equal capacity. Since all edges of $G$ have weight $1$, this proves \textup{(i)}.

For $xy\in E(T)$, let $A$ be the vertex set of either component of $T-xy$. The cut determined by $A$ in $T$ consists only of $xy$ and is a minimum $(x,y)$-cut. Its capacity in $G$ is therefore also $w(xy)$, giving
$$d_G(A)=\sum_{e\in\partial_G(A)}1=w(xy).$$
This proves \textup{(ii)}. Since $G$ is connected, $d_G(A)$ is a positive integer, so all edge weights of $T$ are positive integers.
\end{proof}

\begin{lemma}\label{lem:cut}
Let $d\ge1$ be an integer and $G$ be a connected graph with $\delta(G)\ge d$. There is a set $A\subseteq V(G)$ with $\lvert A\rvert\ge d+1$ such that every cut of $G$ that separates two vertices of $A$ has at least $d$ edges, and $G[A]$ has average degree greater than $d-1$.
\end{lemma}

\begin{proof}
Since $\delta(G)\ge d\ge1$, the graph $G$ has at least two vertices. Choose a weighted tree $T$ as in Lemma~\ref{lem:gomory-hu}.

Contract each tree edge of weight at least $d$. If the resulting tree has exactly one vertex, let $A=V(G)$. Otherwise, choose a leaf $u$ of the resulting tree and let $A$ be the set of vertices of $T$ that were contracted into $u$. For any two distinct vertices of $A$, the path joining them in $T$ consists entirely of edges that were contracted. Each of these edges has weight at least $d$. By Lemma~\ref{lem:gomory-hu}(i), every edge cut of $G$ separating two vertices of $A$ therefore has at least $d$ edges.
If $A\ne V(G)$, then exactly one edge of $T$ joins $A$ to $V(G)\setminus A$, since $A$ was contracted to a leaf. This edge was not contracted, so its weight is less than $d$. Deleting it from $T$ leaves two components with vertex sets $A$ and $V(G)\setminus A$. By Lemma~\ref{lem:gomory-hu}(ii), $d_G(A)$ equals the weight of this edge, and hence $d_G(A)\le d-1$. If $A=V(G)$, then $d_G(A)=0$, so the same inequality holds.

Let $a=|A|$. If $a\le d$, then each vertex of $A$ has at most $a-1$ neighbors in $A$ and hence at least $d-a+1$ neighbors outside $A$. Thus, $d_G(A)\ge a(d-a+1)=d+(a-1)(d-a)\ge d$, contrary to $d_G(A)\le d-1$. Therefore $a\ge d+1$.
Summing the degrees in $G$ over $A$ gives $2e(G[A])+d_G(A)\ge da$. Since $d_G(A)\le d-1$, we obtain $2e(G[A])\ge da-(d-1)$. Dividing by $a$, we know that $G[A]$ has average degree at least $d-(d-1)/a>d-1$, since $a\ge d+1$.
\end{proof}

\begin{proof}[Proof of Theorem~\ref{thm:edge}]
For $m\le1$, the required copy uses no edges and exists since $G$ has at least $k+1$ vertices. Let $d=\max\{k+e,m-1\}$ and choose $A$ by Lemma~\ref{lem:cut}. Join the components of $F$ by adding edges between them to obtain an $m$-vertex tree. Since $G[A]$ has average degree greater than $d-1\ge m-2$, Lemma~\ref{lem:erdos-sos} gives a copy of this tree in $G[A]$. Retain the edges corresponding to the original forest and all $m$ vertices to obtain $F'$.

Let $W$ be a nonempty proper subset of $V(G)$ and consider the cut $\partial_G(W)$. If $A$ is entirely on one side, no edge of $F'$ crosses the cut, so at least $k$ edges remain. If the cut separates two vertices of $A$, it has at least $d$ edges by Lemma~\ref{lem:cut}. Deleting the $e$ edges of $F'$ therefore leaves at least $d-e\ge k$ edges in the cut. Thus $G-E(F')$ is $k$-edge-connected.
\end{proof}

\section*{Acknowledgment}
The second author is grateful to Qinghai Liu for sharing some documents and his knowledge on Mader's conjecture. He is also grateful to Yingzhi Tian from whom he was interested in Mader's conjecture and related topics.

\section*{Declaration on the Use of AI}
The authors used ChatGPT 5.6 to assist in discussing proof strategies, checking proofs, and
improving exposition. In particular, AI recommends the reference \cite{abboud-gh2021}.

\end{document}